\documentclass[preprint,12pt]{article}

\usepackage{amssymb, amsmath, amsthm,graphicx}

\newtheorem{thm}{Theorem}[section]

\newtheorem{cor}[thm]{Corollary}
\newtheorem{lemma}[thm]{Lemma}

\theoremstyle{definition}

\begin{document}

\title{Exit times for multivariate \\ autoregressive processes}
\author{Brita Jung}
\date{}
\maketitle

\begin{center}
Department of Mathematics, \AA bo Akademi University,\\ FIN-20500 \AA bo, Finland
\end{center}

\begin{abstract}
We study exit times from a set for a family of multivariate autoregressive processes with normally distributed noise. By using the large deviation principle, and other methods, we show that the asymptotic behavior of the exit time depends only on the set itself and on the covariance matrix of the stationary distribution of the process. The results are extended to exit times from intervals for the univariate autoregressive process of order $n$, where the exit time is of the same order of magnitude as the exponential of the inverse of the variance of the stationary distribution. 

\end{abstract}


\section{Introduction}

We consider a multivariate autoregressive process with normally distributed noise, defined through
\[
X_t^\varepsilon = AX_{t-1}^\varepsilon+\varepsilon\xi_t,\; t\ge 1,\; X_0^\varepsilon = x_0,
\]
where $X_t^\varepsilon\in \mathbb{R}^d$, $A$ is a real $d\times d$ matrix , $\varepsilon$ is a small positive parameter and $\{\xi_t\}_{t\ge 1}$ is an i.i.d. sequence of multivariate standard normal random variables. We will study the time until the process exits from a set of the type $\{x\in \mathbb{R}^d: |c^T x| < 1\}$ for some vector $c\in \mathbb{R}^d$. Subject to some conditions, we  show that the expectation of this exit time is of the order of magnitude $\exp(1/(\varepsilon^2 c^T\Sigma_\infty  c))$ for small values of $\varepsilon$, where $\varepsilon^2\Sigma_\infty$ is the covariance matrix of the stationary distribution of the process.

\bigskip\noindent
The corresponding univariate case, where $X_t^\varepsilon\in \mathbb{R}\;\forall t\ge 0$ has been investigated before, by Klebaner and Liptser in \cite{Kle} and by Ruths in \cite{Rut}. In \cite{Kle}, the authors proved a large deviation principle (LDP) for a class of past-dependent models. As an example, they used the univariate 
autoregressive process $\{X_t^\varepsilon\}_{t\ge 0}$, where
\[
X_t^\varepsilon = aX_{t-1}^\varepsilon + \varepsilon\xi_t, X_0^\varepsilon = x_0,
\]
where $X_t^\varepsilon\in \mathbb{R}\;\forall t$, $|a|<1$, $\varepsilon$ is a positive parameter and $\{\xi_t\}_{t\ge 1}$ is an i.i.d. sequence of standard normal random variables. This process has a stationary distribution which is normal with mean 0 and variance $1/(1-a^2)$. Klebaner and Liptser showed that the family of processes $\{X_t^\varepsilon\}_{t\ge 0}$ obeys an LDP with rate of speed $\varepsilon^2$ and rate function
\[
I(\bar u) = \left\{\begin{array}{ll}
					\frac12\sum_{t=1}^\infty (u_t-au_{t-1})^2, & u_0 = x_0,\\
					\infty , &\mbox{otherwise,}
				\end{array}\right.
\]
where $\bar u = (u_0,u_1,\ldots )$, and that this implies that
\begin{equation}
\limsup_{\varepsilon\rightarrow 0}\varepsilon^2\log E\tau^\varepsilon \le\frac12 (1-a^2),
\end{equation}
when $\tau^\varepsilon := \min\{t\ge 1: |X_t^\varepsilon|\ge 1\}$. This upper bound is sharp, considering \cite{Rut}, where the corresponding lower bound was proven with another method. We also note the correspondence between this bound and the variance of the stationary distribution. 

\bigskip\noindent

In section 2 of this paper, we establish the corresponding large deviation principle for a family of multivariate processes. We also present a method to get a lower bound for the exit time of normally distributed processes. In section 3 we prove the asymptotics of the exit time of the multivariate autoregressive process. In section 4, we apply the same methods to get a result for the exit time from an interval for the univariate autoregressive process of order $n$, where
\[
X_t^\varepsilon = a_1X_{t-1}^\varepsilon+\ldots + a_nX_{t-n}^\varepsilon + \varepsilon\xi_t,\; t\ge n, \; X_0^\varepsilon = x_0,\ldots , X_{n-1}^\varepsilon = x_{n-1},
\]
where $a_1,\ldots ,a_n$ are real parameters and $\{\xi_t\}_{t\ge n}$ is a sequence of univariate standard normal random variables. 


\section{Methods for upper and lower bounds}

In the first two parts of this section, we consider how to use the large deviation principle to get an upper bound of the asymptotics of an exit time from a set for a process. In the third part of the section, we consider another method for the corresponding lower bound, when the process has a normal distribution.


\subsection{The large deviation principle}

The following definition of the large deviation principle is taken from Varadhan (\cite{Var}), with the slight difference that we let the rate of speed be a function of $\varepsilon$ and call it $q(\varepsilon)$, as Klebaner and Liptser also did (\cite{Kle}). The large deviation principle (LDP) is then defined in the following way: Let $\{P_\varepsilon\}$ be a family of probability measures on the Borel subsets of a complete separable metric space $Z$. We say that $\{P_\varepsilon\}$ satisfies the large deviation principle with a rate function $I(\cdot)$ if there exists a function $I$ from $Z$ into $[0,\infty]$ satisfying the following conditions: $0\le I(z)\le\infty\;\forall z\in Z$, $I$ is lower semicontinuous, the set $\{z:I(z)\le l\}$ is a compact set in  $Z$ for all $\;l<\infty$ and
\begin{eqnarray*}
& &\limsup_{\varepsilon\rightarrow 0} q(\varepsilon )\log P_\varepsilon (C) \le -\inf_{z\in C}I(z) \;\; \mbox{ for every closed set } C\subset Z \mbox{ and }\\
& &\liminf_{\varepsilon\rightarrow 0} q(\varepsilon )\log P_\varepsilon (G) \ge -\inf_{z\in G}I(z) \;\; \mbox{ for every open set } G\subset Z.
\end{eqnarray*}
We will consider a family of processes $\{X_t^\varepsilon\}_{t\ge 0}$, where $X_t^\varepsilon \in \mathbb{R}^d \;\forall t\ge 0$ and
\begin{equation}\label{process_definition}
X_t^\varepsilon = f(X_{t-1}^\varepsilon,\ldots , X_{t-n}^\varepsilon, \varepsilon\xi_t ) \mbox{ for }t\ge n,
\end{equation}
where $f: (\mathbb{R}^d)^{n+1}\mapsto \mathbb{R}^d$ is a continuous function, $\{\xi_t\}_{t\ge n}$ is an i.i.d. sequence of random variables in $\mathbb{R}^d$, $\varepsilon$ is a positive parameter and the starting values are $X_0^\varepsilon =x_0,\ldots , X_{n-1}^\varepsilon= x_{n-1}$. We will prove a large deviation principle for the family of probability measures induced by $\{X_t^\varepsilon\}_{t\ge 0}$, assuming that a large deviation principle for the family of probability measures induced by $\{\varepsilon\xi_n\}$ holds. 

\begin{thm}\label{thm}
Assume that the family of probability measures induced by $\{\varepsilon\xi\}$, where $\xi$ is a copy of $\xi_n$, satisfies a large deviation principle with rate function $I_{\varepsilon\xi}(z)$ and rate of speed $q(\varepsilon)$. Then the large deviation principle holds for the family of probability measures induced by $\{X_t^\varepsilon\}_{t\ge 0}$ with the same rate of speed and the rate function 
\[
I(y_0,y_1,y_2,\ldots ) = \inf_{\genfrac{}{}{0pt}{1}{\genfrac{}{}{0pt}{1}{z_t\in \mathbb{R}^d\;\forall t\ge n}{y_t = f(y_{t-1},\ldots , y_{t-n},z_t), t\ge n}}{y_0 = x_0,\ldots , y_{n-1}= x_{n-1}}} \sum_{t=n}^\infty I_{\varepsilon\xi}(z_t).
\]
\end{thm}
\noindent
Proof: We have assumed that an LDP holds for the family of probability measures induced by the family $\{\varepsilon\xi\}$, with the rate function $I_{\varepsilon\xi}(z)$ and the rate of speed $q(\varepsilon)$. By \cite{Lyn}, the LDP then holds for the family of probability measures induced by the family of vectors $\{\varepsilon\xi_t\}_{t=n}^N$ with the same rate of speed and the rate function
\[
I_{\{\varepsilon\xi_t\}_{t=n}^N}(z_n,\ldots , z_N) = \sum_{t=n}^{N} I_{\varepsilon\xi}(z_t),
\]
where $N$ is finite. By the Dawson-G\"{a}rtner theorem (see for example \cite{Dem}), it follows that the LDP holds for the family of probability measures induced by $\{\varepsilon\xi_t\}_{t\ge n}$ with rate of speed $q(\varepsilon)$ and rate function
\[
I_{\{\varepsilon\xi_t\}_{t\ge n}}(z_n,z_{n+1},\ldots) = \sum_{t=n}^\infty I_{\varepsilon\xi_t}(z_t).
\]
Now, since $f$ is continuous, the mapping $\{\varepsilon\xi_t\}_{t\ge n}\mapsto \{X_t^\varepsilon\}_{t\ge 0}$ is continuous in the space $(\mathbb{R}^d)^\infty$ with the metric $\rho(x,y) = \sum_{j\ge 1}2^{-j}\frac{||x_j-y_j||}{1+ ||x_j-y_j||}$, where $||\cdot ||$ denotes the Euclidian norm on $\mathbb{R}^d$. Thus, we can use the contraction principle (see for example \cite{Dem}). It implies that the LDP for the family of probability measures associated with the family $\{X_t^\varepsilon\}_{t\ge 0}$ holds with rate of speed $q(\varepsilon)$ and rate function
\[
I(y_0, y_1,y_2,\ldots ) = \inf_{\genfrac{}{}{0pt}{1}{\genfrac{}{}{0pt}{1}{(z_n,z_{n+1},\ldots )\in (\mathbb{R}^d)^\infty }{y_t = f(y_{t-1},\ldots , y_{t-n},z_t),t\ge n}}{y_0 = x_0,\ldots , y_{n-1}= x_{n-1}}}I_{ \{ \varepsilon\xi_t \}_{t\ge n}}(z_n,z_{n+1},\ldots ),
\]
where the infimum over the empty set is taken as $\infty$. We can write this rate function as
\[
I(y_0,y_1,y_2,\ldots ) = \inf_{\genfrac{}{}{0pt}{1}{\genfrac{}{}{0pt}{1}{z_t\in \mathbb{R}^d\;\forall t\ge n}{y_t = f(y_{t-1},\ldots , y_{t-n},z_t),t\ge n}}{y_0 = x_0,\ldots , y_{n-1}= x_{n-1}}} \sum_{t=n}^\infty I_{\varepsilon\xi}(z_t)
\]
and the proof is finished.


\subsection{Exit times with the large deviation principle}

We will later use the large deviation principle to get a bound for the exit time from a set for a certain process. To see how this will be done, let us for the moment define the exit time as
\begin{equation}
\tau := \min\{t\ge n: X_t^\varepsilon \notin \Omega\},
\end{equation}
where $X_t^\varepsilon$ is defined as in equation \ref{process_definition} and $\Omega$ is a set in $\mathbb{R}^d$. Assume that the starting points $x_0,\ldots ,x_{n-1}$ of the process belong to $\Omega$. For the expectation of the exit time, we have the following, where $M$ is any integer greater than or equal to $n$:
\begin{eqnarray*}
E_{x_0,\ldots ,x_{n-1}}(\tau )&\le & M + P(\tau > M-1)E_{x_0,\ldots ,x_{n-1}}(\tau | \tau > M-1)\\
&\le & M + P(\tau >M-1)[MP(\tau = M|\tau > M-1) \\
& & + (M+\sup_{x_0,\ldots ,x_{n-1}\in \Omega}E_{x_0,\ldots ,x_{n-1}}(\tau ))P(\tau > M|\tau > M-1)]\\
&\le & 2M + \sup_{x_0,\ldots ,x_{n-1}\in \Omega} E_{x_0,\ldots ,x_{n-1}}(\tau )\cdot \sup_{x_0,\ldots ,x_{n-1}\in \Omega}P_{x_0,\ldots ,x_{n-1}}(\tau > M).
\end{eqnarray*} 
Thus, it holds that
\[
\sup_{x_0,\ldots ,x_{n-1}\in \Omega} E_{x_0,\ldots ,x_{n-1}}(\tau ) \le \frac{2M}{\inf_{x_0,\ldots ,x_{n-1}\in \Omega}P_{x_0,\ldots ,x_{n-1}}(\tau\le M)},
\]
or simply that
\begin{equation}
E_{x_0,\ldots ,x_{n-1}}(\tau ) \le \frac{2M}{\inf_{x_0,\ldots ,x_{n-1}\in \Omega}P_{x_0,\ldots ,x_{n-1}}(\tau\le M)}
\end{equation}
for any set of starting points $x_0,\ldots ,x_{n-1}\in \Omega$ and any integer $M\ge n$. If the infimum in the denominator is attained for the starting points $x_0^*,\ldots ,x_{n-1}^*\in\Omega$, the inequality above implies that
\begin{equation}\label{inequality}
\limsup_{\varepsilon\rightarrow 0} q(\varepsilon )\log E_{x_0,\ldots ,x_{n-1}}(\tau) \le -\lim_{\varepsilon\rightarrow 0} q(\varepsilon )\log P_{x_0^*,\ldots ,x_{n-1}^*}(\tau\le M),
\end{equation}
if the right hand side limit exists. Since
\[
P_{x_0^*,\ldots ,x_{n-1}^*}(\tau\le M) = P_{x_0^*,\ldots ,x_{n-1}^*}(X_t^\varepsilon \notin \Omega \mbox{ for some } t\in \{n,\ldots , M\}),
\]
the limit may be calculated if we have a large deviation principle for the family of probability measures induced by $\{X_t^\varepsilon\}_{t\ge 0}$ and if the function $f$ and the set $\Omega$ are suitable. 

In sections 3 and 4 we will use this method to get upper bounds for exit times for multivariate autoregressive processes and univariate processes of order $n$, respectively.


\subsection{A lower bound for the exit time of normally distributed variables}

We now leave the large deviation principle for a moment, and consider a method to get  a lower bound for an exit time. The following theorem gives a lower bound for the asymptotics of the mean exit time from a symmetric interval for a sequence of univariate normally distributed random variables $\{Y_t^\varepsilon\}_{t\ge 1}$, with mean zero and bounded variance. Thus, in this section we consider the exit time
\[
\tau_{(-1,1)} := \min\{t\ge 1: |Y_t^\varepsilon |\ge 1\}.
\]

\begin{thm}\label{thmlower}
Assume that $\{Y_t^\varepsilon\}_{t\ge 1}$ is a sequence of normally distributed random variables, all with mean 0, and that
\[
{\rm Var}(Y_t^\varepsilon ) \le q(\varepsilon )\sigma^2 \;\forall t\ge 1,
\]
for some $\sigma^2 > 0$ and some positive function $q(\varepsilon )$ where $\lim_{\varepsilon\rightarrow 0} q(\varepsilon) = 0$. Then 
\[
\liminf_{\varepsilon\rightarrow 0}q(\varepsilon )\log E\tau_{(-1,1)} \ge \frac{1}{2\sigma^2}.
\]
\end{thm}
\noindent
Proof: Since $Y_t^\varepsilon$ has a normal distribution with mean zero and variance bounded by $q(\varepsilon )\sigma^2$, $E(e^{\lambda Y_t^\varepsilon})\le \exp(\frac12 \lambda^2q(\varepsilon )\sigma^2)$  $\forall t\ge 1$, which implies that also
\[
E(\cosh(\lambda Y_t^\varepsilon)) \le e^{\frac12 \lambda^2q(\varepsilon )\sigma^2}\;\; \forall t\ge 1.
\]
For any $N\ge 1$, we have the following Chernoff-type bound of the probability that the exit time is smaller than or equal to $N$:
\begin{eqnarray*}
 P(\tau_{(-1,1)} \le N) &=& P(\max_{1\le t\le [N]}|Y_t^\varepsilon |\ge 1)
= P(\cosh (\lambda\max_{1\le t\le [N]}|Y_t^\varepsilon | )\ge \cosh\lambda )\\
& & \le (\cosh\lambda)^{-1} E(\cosh (\lambda\max_{1\le t\le [N]}|Y_t^\varepsilon | )),
\end{eqnarray*}
which holds for any positive $\lambda$. (Of course, the bound holds for any $\lambda\in R$.) Since
\[
\cosh (\lambda\max_{1\le t\le [N]}|Y_t^\varepsilon | ) = \max_{1\le t\le [N]}\cosh(\lambda Y_t^\varepsilon) \le \sum_{t=1}^{[N]}\cosh(\lambda Y_t^\varepsilon),
\] 
it follows that
\[
E(\cosh (\lambda\max_{1\le t\le [N]}|Y_t^\varepsilon | )) \le [N]e^{\frac12 \lambda^2q(\varepsilon )\sigma^2} \le Ne^{\frac12 \lambda^2q(\varepsilon )\sigma^2}.
\]
Thus, we have the bound
\[
P(\tau_{(-1,1)}\le N) \le (\cosh\lambda)^{-1} Ne^{\frac12 \lambda^2q(\varepsilon )\sigma^2} \le 2e^{-\lambda}Ne^{\frac12 \lambda^2q(\varepsilon )\sigma^2},
\]
for any $\lambda > 0$. By choosing $\lambda$ in the optimal way, that is, as $\lambda = 1/(q(\varepsilon )\sigma^2)$, we get the bound
\[
P(\tau_{(-1,1)}\le N) \le 2N \exp(-\frac{1}{2q(\varepsilon ) \sigma^2}).
\]
Now, let $\delta$ be a small positive number and choose $N = \exp(\frac{1}{2q(\varepsilon ) \sigma^2}-\frac{\delta}{q(\varepsilon)})$. Then $P(\tau_{(-1,1)} >N) > 1-2\exp(-\frac{\delta}{q(\varepsilon )})$, which implies that
\[
E\tau_{(-1,1)} \ge NP(\tau_{(-1,1)} > N) \ge \exp(\frac{1}{2q(\varepsilon ) \sigma^2}-\frac{\delta}{q(\varepsilon )})(1-2\exp(-\frac{\delta}{q(\varepsilon )})),
\]
and thus
\[
\liminf_{\varepsilon\rightarrow 0}q(\varepsilon )\log E\tau_{(-1,1)} \ge \frac{1}{2\sigma^2}-\delta.
\]
Since this holds for any $\delta >0$, we get the lower bound
\[
\liminf_{\varepsilon\rightarrow 0}q(\varepsilon )\log E\tau_{(-1,1)} \ge \frac{1}{2\sigma^2},
\]
and the proof is finished.

\medskip
\noindent
{\bf Remark:} If we wanted to consider a one-sided exit time, for example the time until $Y_t^\varepsilon > 1$, we could simply use the exponential function instead of the hyperbolic cosine in the argument above. The resulting lower bound would be the same.


\section{Exit times for a multivariate autoregressive process}

In this section, we use the methods described in section 2 to show that the expectation of the exit time from a set $\{x\in \mathbb{R}^d: |c^T x| < 1\}$ for a vector $c\in \mathbb{R}^d$ (that is not the zero vector) for a multivariate autoregressive process is of the order of magnitude $\exp(1/(\varepsilon^2 c^T\Sigma_\infty  c))$, where $\varepsilon^2\Sigma_\infty$ is the covariance matrix of the stationary distribution of the process.


\subsection{A multivariate autoregressive process}

By a multivariate autoregressive process, we mean a process $\{X_t^\varepsilon\}_{t\ge 0}$, such that
\begin{equation}\label{multivariate_ar}
X_t^\varepsilon = AX_{t-1}^\varepsilon + \varepsilon\xi_t, \; X_0^\varepsilon = x_0,
\end{equation}
where $X_t^\varepsilon \in \mathbb{R}^d\;\forall t$, $A$ is a real $d\times d$ matrix, $\varepsilon$ is a positive parameter and $\{\xi_t\}_{t\ge 1}$ is an i.i.d. sequence of multivariate normal random variables in $\mathbb{R}^d$, with mean zero and covariance matrix $I$ (the unit matrix). For any $t\ge 1$, $X_t^\varepsilon$ has a multivariate normal distribution with mean $EX_t^\varepsilon = A^tx_0$, where $x_0$ is the starting point,  and covariance matrix $\varepsilon^2\Sigma_t$, where
\begin{equation}
\Sigma_t = A\Sigma_{t-1}A^T + I, \;\; t\ge 2, 
\end{equation}
and $\Sigma_1 = I$. The matrix $\Sigma_t$ can also be written as the sum
\begin{equation}
\Sigma_t = \sum_{i=0}^{t-1} A^i(A^T)^i, t\ge 1.
\end{equation}
Throughout, we will assume that all eigenvalues of $A$ have absolute values less than one. The process $\{X_t^\varepsilon\}_{t\ge 0}$ then has a stationary distribution,  which is multivariate normal with mean $(0,0,\ldots ,0)^T$ and covariance matrix $\varepsilon^2\Sigma_\infty$, where $\Sigma_\infty$ satisfies
\begin{equation}\label{covmatrixequation}
\Sigma_\infty =  A\Sigma_\infty A^T + I. 
\end{equation}
Of course, the matrix $\Sigma_\infty$ can also be expressed as the sum
\begin{equation}
\Sigma_\infty = \sum_{i=0}^\infty A^i(A^T)^i.
\end{equation}


\subsection{Exit times for the multivariate autoregressive process}

For the multivariate autoregressive process $\{X_t^\varepsilon\}_{t\ge 0}$, we will consider the exit time
\begin{equation}
\tau:=\min\{t\ge 1: |c^TX_t^\varepsilon |\ge 1\},
\end{equation}
where $c$ is a vector in $\mathbb{R}^d$, $c\neq (0,\ldots , 0)^T$. We will find the limit of $\varepsilon^2\log E\tau$ as $\varepsilon\rightarrow 0$, by using the methods described in section 2. For the upper bound, we will use the large deviation principle, so we need the following corollary.

\begin{cor}\label{corollary}
The family of probability measures induced by $\{X_t^\varepsilon\}_{t\ge 0}$, where $X_t^\varepsilon$ is defined as in equation \ref{multivariate_ar}, satisfies the large deviation principle  with rate of speed $q(\varepsilon) = \varepsilon^2$ and rate function
\[
I(y_0,y_1,\ldots ) = \frac12 \sum_{t=1}^\infty (y_t-Ay_{t-1})^T(y_t-Ay_{t-1}),
\]
where $y_0 = x_0$.
\end{cor}
\noindent
Proof: By using Cram\'{e}r's theorem (see for example \cite{Dem}), one can show that the family of probability measures induced by the family $\{\varepsilon\xi\}$, where $\xi$ is multivariate normal with mean zero and covariance matrix $I$, satisfies the LDP with rate of speed $\varepsilon^2$ and rate function
\[
I_{\varepsilon\xi} (z) = \frac12 z^Tz, \; z\in \mathbb{R}^d.
\]
By using Theorem \ref{thm}, we can deduce that the family of probability measures induced by $\{X_t^\varepsilon\}_{t\ge 0}$ satisfies the LDP with the same rate of speed $\varepsilon^2$ and rate function
\begin{eqnarray*}
I(y_0,y_1,\ldots ) &=& \inf_{\genfrac{}{}{0pt}{1}{\genfrac{}{}{0pt}{1}{z_t\in \mathbb{R}^d}{y_t = Ay_{t-1} + z_t}}{y_0 = x_0}} \sum_{t=1}^\infty \frac12 z_t^Tz_t \\
&=& \frac12\sum_{t=1}^\infty (y_t-Ay_{t-1})^T(y_t-Ay_{t-1}),
\end{eqnarray*}
where $y_0 = x_0$.

\medskip\noindent
For the exit time of a multivariate autoregressive process, starting at the origin, we will prove the following theorem:
\begin{thm}\label{ar_thm}
For the exit time $\tau = \min\{t\ge 1: |c^TX_t^\varepsilon|\ge 1\}$, where $\{X_t^\varepsilon\}_{t\ge 0}$ is the multivariate autoregressive process defined in equation \ref{multivariate_ar}, and $x_0 = (0,\ldots ,0)^T$,
\[
\lim_{\varepsilon\rightarrow 0}\varepsilon^2\log E\tau = \frac{1}{2c^T\Sigma_\infty c},
\]
where $\varepsilon^2\Sigma_\infty$ is the covariance matrix of the stationary distribution of the process.
\end{thm}

\noindent
Proof: The theorem follows from lemmas \ref{upper} and \ref{lower} below.

\begin{lemma}\label{upper}
For $\{X_t^\varepsilon\}_{t\ge 0}$ and $\tau$ as in Theorem \ref{ar_thm}, we have
\[
\limsup_{\varepsilon\rightarrow 0}\varepsilon^2\log E_{x_0}\tau \le \frac{1}{2c^T\Sigma_\infty c},
\]
for any $x_0$ such that $|c^Tx_0|<1$.
\end{lemma}
\noindent
Proof: Consider the exit time $\tau = \min\{t\ge 1: |c^TX_t^\varepsilon|\ge 1\}$. This means that we consider exits from the set $\Omega := \{x\in \mathbb{R}^d: |c^Tx| < 1\}$. For this set $\Omega$, $\inf_{x_0\in \Omega}P_{x_0}(\tau\le M) = P_{(0,\ldots ,0)^T}(\tau\le M)$, since $X_t^\varepsilon$ has a normal distribution with mean $A^tx_0$. Thus, inequality \ref{inequality} in section 2.2 implies that
\[
\limsup_{\varepsilon\rightarrow 0}\varepsilon^2 \log E_{x_0}\tau \le -\lim_{\varepsilon\rightarrow 0}\varepsilon^2\log P_{(0,\ldots ,0)^T}(\tau\le M),
\]
where the right hand side limit can be calculated with the LDP that was proven in corollary \ref{corollary}. Since
\[
\{\tau\le M\} = \{\max_{1\le t\le M}|c^TX_t^\varepsilon |\ge 1\},
\]
we have
\[
\lim_{\varepsilon\rightarrow 0}\varepsilon^2\log P_{(0,\ldots ,0)^T}(\tau\le M) = -\inf_{\genfrac{}{}{0pt}{1}{\max_{1\le t\le M}|c^Ty_t |\ge 1,}{y_0 = (0,\ldots ,0)^T}} \frac12\sum_{t=1}^\infty (y_t-Ay_{t-1})^T(y_t-Ay_{t-1}).
\]
Consider this infimum. The following holds:
\begin{eqnarray*}
& &\inf_{\genfrac{}{}{0pt}{1}{\max_{1\le t\le M}|c^Ty_t |\ge 1}{y_0 = (0,\ldots ,0)^T}} \frac12\sum_{t=1}^\infty (y_t-Ay_{t-1})^T(y_t-Ay_{t-1})\\
&=& \inf_{1\le N\le M}\left( \inf_{\genfrac{}{}{0pt}{1}{|c^Ty_N |\ge 1}{y_0 = (0,\ldots ,0)^T}} \frac12\sum_{t=1}^\infty (y_t-Ay_{t-1})^T(y_t-Ay_{t-1})\right)\\
&=& \inf_{1\le N\le M}\left( \inf_{\genfrac{}{}{0pt}{1}{|c^Ty_N |\ge 1}{y_0 = (0,\ldots ,0)^T}} \frac12\sum_{t=1}^N (y_t-Ay_{t-1})^T(y_t-Ay_{t-1})\right),
\end{eqnarray*}
where the last equality holds because we can choose $y_t = Ay_{t-1}$ for $t>N$. We can write $y_N$ as the telescoping sum $\sum_{t=1}^N A^{N-t}(y_t-Ay_{t-1})$, when $y_0=(0,\ldots ,0)^T$. By using the Cauchy-Schwarz inequality, we get
\begin{eqnarray*}
& & \left(\sum_{t=1}^N (y_t-Ay_{t-1})^T(y_t-Ay_{t-1})\right) \cdot \left(\sum_{t=1}^N c^TA^{N-t}(A^{N-t})^Tc \right)\\
&\ge & \left(\sum_{t=1}^N c^TA^{N-t}(y_t-Ay_{t-1})\right)^2 = (c^Ty_N)^2.
\end{eqnarray*}
Equality in the Cauchy-Schwarz inequality is attained when $y_t-Ay_{t-1} = K(A^{N-t})^Tc$, $t=1,\ldots , N$, for any constant $K\in \mathbb{R}$. This holds when
\[
y_t = K\left(\sum_{i=0}^{t-1}A^i(A^T)^i\right)(A^{N-t})^Tc = K\Sigma_t(A^{N-t})^Tc, \mbox{  for  } t= 1,\ldots , N,
\]
where $\Sigma_t$ is defined as in section 3.1. By choosing $K = 1/(|c^T\Sigma_Nc|)$, we get $|c^Ty_N|= 1$. Thus, we have now shown that
\begin{eqnarray*}
\inf_{\genfrac{}{}{0pt}{1}{|c^Ty_N |\ge 1,}{y_0 = (0,\ldots ,0)^T}} \frac12\sum_{t=1}^N (y_t-Ay_{t-1})^T(y_t-Ay_{t-1})
= \frac{1}{2\sum_{t=1}^N c^TA^{N-t}(A^T)^{N-t}c} = \frac{1}{2c^T\Sigma_N c}.
\end{eqnarray*}
Since $\Sigma_t = \sum_{i=0}^{t-1}A^i(A^T)^i$ and
\[
c^T\Sigma_tc = c^T\Sigma_{t-1}c + c^TA^{t-1}(A^{t-1})^Tc \ge c^T\Sigma_{t-1}c \;\;\forall t = 2,\ldots , M,
\]
$\{c^T\Sigma_t c\}_{t\ge 1}$ is a positive and increasing sequence. It follows, that
\[
\inf_{1\le N\le M} \frac{1}{2c^T\Sigma_Nc} = \frac{1}{2c^T\Sigma_Mc},
\]
and we have shown that
\[
\limsup_{\varepsilon\rightarrow 0}\varepsilon^2\log E_{x_0}\tau \le \frac{1}{2c^T\Sigma_Mc}.
\]
Since this inequality holds for any integer $M\ge 1$, we actually have
\begin{equation}
\limsup_{\varepsilon\rightarrow 0}\varepsilon^2\log E_{x_0}\tau \le \frac{1}{2c^T\Sigma_\infty c},
\end{equation}
and the proof is finished.

\begin{lemma}\label{lower}
For $\{X_t^\varepsilon\}_{t\ge 0}$ and $\tau$ as in Theorem \ref{ar_thm}, we have
\[
\liminf_{\varepsilon\rightarrow 0}\varepsilon^2\log E_{(0,\ldots ,0)^T}(\tau ) \ge \frac{1}{2c^T\Sigma_\infty c},
\]
where $E_{(0\ldots ,0)^T}$ denotes that the starting point of the process is $x_0 = (0,\ldots ,0)^T$.
\end{lemma}

\noindent
Proof: For each $t\ge 1$, $c^TX_t^\varepsilon$ has a univariate normal distribution with mean zero (since $x_0$ is now chosen to be the zero vector) and variance ${\rm Var}(c^TX_t^\varepsilon) = \varepsilon^2c^T\Sigma_t c$. In the proof of lemma \ref{upper}, we showed that $\{c^T\Sigma_t c\}_{t\ge 1}$ is an increasing sequence, and thus we have
\[
{\rm Var}(c^TX_t^\varepsilon) \le \varepsilon^2c^T\Sigma_\infty c \;\;\forall t\ge 1.
\]
The statement of the lemma then follows immediately from theorem \ref{thmlower}.


\bigskip\noindent
We illustrate the result in theorem \ref{ar_thm} by simulating a bivariate process $\{X_t^\varepsilon\}_{t\ge 0}$, where
\[
X_t^\varepsilon = AX_{t-1}^\varepsilon+\varepsilon\xi_t,
\]
where $X_t^\varepsilon = (X_{t,1}^\varepsilon,X_{t,2}^\varepsilon)^T\in \mathbb{R}^2 \;\forall t\ge 1$, $X_0 = (0,0)^T$, $\{\xi_t\}_{t\ge 1}$ is an i.i.d. sequence of bivariate standard normal random variables and $A = \left( \begin{array}{cc}0.8 & 1\\ 0 & 0.5\end{array}\right)$. Since the eigenvalues of $A$ (0.8 and 0.5) are less than one in absolute value, the process has a stationary distribution. Let $c=(1,1)^T$ and consider the exit time $\tau = \min\{t\ge 1: |X_{t,1}^\varepsilon+X_{t,2}^\varepsilon|\ge 1\}$. The matrix $\Sigma_\infty$ is calculated from equality \ref{covmatrixequation}. We get
\[
\Sigma_\infty = \left( \begin{array}{cc}\frac{925}{81} & \frac{10}{9}\\ \frac{10}{9} & \frac{12}{9}\end{array}\right). 
\]
Theorem \ref{ar_thm} says that $\lim_{\varepsilon\rightarrow 0}\varepsilon^2\log E\tau = 1/(2c^T\Sigma_\infty c) = 81/2426 \approx 0.03339$. We use the statistical programming package R to simulate paths of the process for a few values of $\varepsilon$. For each value of $\varepsilon$, 100 paths are simulated and the mean exit time is calculated. The results are shown in table \ref{thetable}. For $\varepsilon = 0.12$, the mean exit time is around 84, while it is around 6 000 000 for $\varepsilon = 0.05$. Naturally, the simulations become more and more time-consuming as $\varepsilon$ decreases and the mean exit time increases.

\begin{table}[ht]
\begin{center}
  \begin{tabular}{ |l|c|c|c|c|c|c|c|}
    \hline
    $\varepsilon$  & 0.12 & 0.10  & 0.08 & 0.07 & 0.06 & 0.05 \\ 
	\hline
    $\varepsilon^2\log E\tau$  & 0.0639 & 0.0554  & 0.0473 & 0.0434 & 0.0415 & 0.0389 \\ 
    \hline
  \end{tabular}
  \caption{Simulation of a bivariate autoregressive process}
  \label{thetable}
\end{center}
\end{table}


\section{Exit times for the autoregressive process of order $n$}

We will now use the methods in section 2 for the univariate autoregressive process of order $n$ with normally distributed noise.

\subsection{The autoregressive process of order $n$}

The autoregressive process of order $n$ is defined as the process $\{X_t^\varepsilon\}_{t\ge 0}$, where
\begin{equation}
X_t^\varepsilon = b_1X_{t-1}^\varepsilon+\ldots + b_nX_{t-n}^\varepsilon + \varepsilon\xi_t,\; X_0^\varepsilon = x_0, \ldots , X_{n-1}^\varepsilon = x_{n-1}.
\end{equation}
Here $X_t^\varepsilon\in \mathbb{R}\;\forall t\ge 0$, $b_1,\ldots , b_n$ are real parameters,  $\varepsilon$ is a positive parameter and $\{\xi_t\}_{t\ge n}$ is an i.i.d. sequence of standard normal (univariate) random variables. We consider the exit time from the interval $(-1,1)$, that is,
\begin{equation}\label{ar_n_tau}
\tau_{(-1,1)} = \min\{t\ge n: |X_t^\varepsilon|\ge 1\}.
\end{equation}
The process can actually be seen as a multivariate process. Let $Y_t^\varepsilon := (X_t^\varepsilon,\ldots , X_{t-n+1}^\varepsilon)^T$ $\forall t\ge n-1$. Then $\{Y_t^\varepsilon\}_{t\ge n}$ is a multivariate process that satisfies
\[
Y_t^\varepsilon = BY_{t-1}^\varepsilon + \varepsilon (\xi_t,0,\ldots ,0)^T\;\mbox{ for } t\ge n, \; Y_{n-1}^\varepsilon = (x_{n-1},\ldots ,x_0)^T,
\]
where 
\[
B = \left(\begin{array}{cccc} b_1 & b_2 &\cdots & b_n\\
								1 & 0 &\cdots &0\\
								0 & \ddots & 0 & 0\\
								0 &\cdots & 1 & 0\end{array}\right).
\]
This process is similar to but not exactly like the multivariate autoregressive process that we considered in section 3. For each $t\ge n$, $Y_t^\varepsilon$ has a multivariate normal distribution with mean $EY_t^\varepsilon = B^{t-n+1}(x_{n-1},\ldots ,x_0)^T$ and covariance matrix $\varepsilon^2\Sigma_t$, where $\Sigma_t$ is given by
\begin{eqnarray*}
\Sigma_t &=& B\Sigma_{t-1}B^T  + (1,0,\ldots ,0)^T(1,0,\ldots ,0) \;\mbox{ for } t\ge n+1,\\ \Sigma_n &=& (1,0,\ldots ,0)^T(1,0,\ldots ,0),
\end{eqnarray*}
or by the sum
\[
\Sigma_t = \sum_{k=0}^{t-n}B^k(1,0,\ldots ,0)^T(1,0,\ldots ,0) (B^T)^k.
\]
Throughout this section, we make the assumption that $b_1,\ldots ,b_n$ are such that all eigenvalues of the matrix $B$ have absolute values less than one. Then the process $\{Y_t^\varepsilon\}_{t\ge n-1}$ has a stationary distribution which is normal with the zero vector as mean and covariance matrix $\varepsilon^2\Sigma_\infty$, where
\begin{eqnarray*}
\Sigma_\infty &=& B\Sigma_\infty B^T  + (1,0,\ldots ,0)^T(1,0,\ldots ,0), \mbox{ or }\\
\Sigma_\infty &=& \sum_{k=0}^\infty B^k(1,0,\ldots ,0)^T(1,0,\ldots ,0) (B^T)^k.
\end{eqnarray*}
This implies that the original univariate process $\{X_t^\varepsilon\}_{t\ge 0}$ has a stationary distribution which is normal with mean zero and variance $\varepsilon^2\sigma^2$, where $\sigma^2$ is given by
\begin{equation}
\sigma^2 = \sum_{k=0}^\infty (B_{11}^k)^2,
\end{equation}
where $B_{11}^k$ denotes the element at the first row and the first column of the matrix $B^k$.

\subsection{Exit times from an interval}

Under the assumption that the starting points $x_0,\ldots ,x_{n-1}$ are zeroes, we have the following result for the exit time $\tau_{(-1,1)}$:

\begin{thm}
For the autoregressive process $\{X_t^\varepsilon\}_{t\ge 0}$ of order $n$, and the exit time $\tau_{(-1,1)}$, 
\[
\lim_{\varepsilon\rightarrow 0} \varepsilon^2 \log E_{(0,\ldots ,0)}\tau_{(-1,1)} = \frac{1}{2\sigma^2},
\]
assuming that all eigenvalues of $B$ are less than one in absolute value, and that $x_0 = \ldots = x_{n-1}= 0$.
\end{thm}

\noindent
Proof: We use the large deviation principle to get an upper bound of the limit. The logarithmic moment generating function of $(\xi,0,\ldots ,0)^T$, where $\xi$ is a standard normal random variable, is
\[
\Lambda(\lambda) = \log E(e^{\lambda_1\xi}) = \frac{\lambda_1^2}{2},
\]
where $\lambda = (\lambda_1,\ldots ,\lambda_n)^T$. Thus, the Fenchel-Legendre transform of $\Lambda (\lambda)$ is 
\[
\Lambda^*(z) = \sup_{\lambda\in R^n}(\lambda^Tz-\Lambda(\lambda)) = \left\{\begin{array}{l} \frac{z_1^2}{2}, \mbox{ if } z_2 = \ldots = z_n = 0\\
 \infty, \mbox{ otherwise.}\end{array} \right. ,
\]
where $z = (z_1,\ldots ,z_n)^T\in \mathbb{R}^n$. Cram\'{e}r's theorem (see for example \cite{Dem}) now gives us that the family of probability measures induced by the family $\{\varepsilon(\xi,0,\ldots ,0)^T\}$ satisfies the large deviation principle with rate of speed $\varepsilon^2$ and rate function $I_{\varepsilon(\xi,0,\ldots ,0)^T}(z) = \Lambda^*(z), z\in \mathbb{R}^n$. By theorem \ref{thm}, the large deviation principle then holds for the family of probability measures induced by $\{Y_t^\varepsilon\}_{t\ge n-1}$ with rate of speed $\varepsilon^2$ and rate function
\begin{eqnarray*}
I(y_{n-1},y_n,\ldots ) &=& \inf_{\genfrac{}{}{0pt}{1}{\genfrac{}{}{0pt}{1}{z_t\in \mathbb{R}^n}{y_t = By_{t-1}+z_t,\; t\ge n}}{y_{n-1}= (x_{n-1},\ldots ,x_0)^T}}\sum_{t=n}^\infty I_{\varepsilon(\xi,0,\ldots ,0)^T}(z_t)\\
&=& \left\{ \begin{array}{ll}\frac12 \sum_{t=n}^\infty ((y_t-By_{t-1})_1)^2 & \mbox{if } (y_t-By_{t-1})_k = 0\\
 & \forall k=2,\ldots ,n,\; \forall t \ge n\\
& \mbox{and } y_{n-1} = (x_{n-1},\ldots ,x_0)^T\\
 \infty ,& \mbox{otherwise,}\end{array}\right.
\end{eqnarray*}
where $y_{n-1}, y_n,\ldots \in \mathbb{R}^n$ and $(y_t-By_{t-1})_k$ denotes the $k$:th element of the vector $y_t-By_{t-1}$. We now proceed as in the multivariate autoregressive case. This time, we consider exits from the set $\Omega = \{x\in \mathbb{R}^d: |c^Tx|<1\}$ for the vector $c = (1,0,\ldots , 0)^T$. For this $\Omega$, $\inf_{y_{n-1}\in\Omega}P_{y_{n-1}}(\tau_{(-1,1)}\le M) = P_{(0,\ldots ,0)^T}(\tau_{(-1,1)}\le M)$ and
\begin{eqnarray*}
& & \lim_{\varepsilon\rightarrow 0}\varepsilon^2\log P_{(0,\ldots ,0)^T}(\tau_{(-1,1)}\le M) = -\inf_{\genfrac{}{}{0pt}{1}{\max_{n\le t\le M}|c^Ty_t|\ge 1}{y_{n-1} = (0,\ldots 0)^T}} I(y_{n-1},y_n,\ldots ) \\
&=& -\inf_{n\le N\le M}(\inf_{\genfrac{}{}{0pt}{1}{|c^Ty_N|\ge 1}{y_{n-1} = (0,\ldots ,0)^T}}I(y_{n-1},y_n,\ldots )),
\end{eqnarray*}
where
\begin{equation}\label{ar_n_inf}
\inf_{\genfrac{}{}{0pt}{1}{|c^Ty_N|\ge 1}{y_{n-1} = (0,\ldots ,0)^T}} I(y_{n-1},y_n,\ldots ) = \!\!\!\!\!\!\!\!\!\!\!\!\!\!\! \inf_{\genfrac{}{}{0pt}{1}{\genfrac{}{}{0pt}{1}{|c^Ty_N|\ge 1}{y_{n-1} = (0,\ldots ,0)^T}}{(y_t-By_{t-1})_k = 0, \;2\le k\le n, \forall t\ge n}}\!\!\!\!\!\!\!\!\!\!\!\!\!\!\!\frac12 \sum_{t=n}^N (y_t-By_{t-1})^T(y_t-By_{t-1}).
\end{equation}
As in the multivariate autoregressive case, one can use the Cauchy-Schwarz inequality to show that the sum on the right hand side in equality \ref{ar_n_inf} is larger than or equal to $1/(2c^T\Sigma_N c)$, where $\Sigma_N$ is defined as in section 4.1. Equality is achieved for $y_t = K\Sigma_t(B^{N-t})^Tc$, $t=n,\ldots , N$ where the choice of $K = 1/(|c^T\Sigma_Nc|)$ gives $|c^Ty_N| = 1$. It is easy to check that this sequence $\{y_t\}_{t=n}^N$ indeed satisfies $(y_t-By_{t-1})_k = 0$ for $k=2,\ldots , n$ and $t=n,\ldots N$. Thus, we get
\[
\inf_{\genfrac{}{}{0pt}{1}{\genfrac{}{}{0pt}{1}{|c^Ty_N|\ge 1}{y_{n-1} = (0,\ldots ,0)^T}}{(y_t-By_{t-1})_k = 0,\; 2\le k\le n, \forall t\ge n}}\frac12 \sum_{t=n}^N (y_t-By_{t-1})^T(y_t-By_{t-1}) = \frac{1}{2c^T\Sigma_N c},
\]
which implies that
\[
\lim_{\varepsilon\rightarrow 0}\varepsilon^2\log P_{(0,\ldots ,0)^T}(\tau_{(-1,1)}\le M) = -\frac{1}{2c^T\Sigma_M c}.
\]
Inequality \ref{inequality} now gives us that $\limsup_{\varepsilon\rightarrow 0}\varepsilon^2\log E_{(x_0,\ldots ,x_{n-1})}\tau_{(-1,1)}\le 1/(2c^T\Sigma_M c)$, and since this holds for any $M\ge n$, we may substitute $\Sigma_\infty$ for $\Sigma_M.$ Since $c = (1,0,\ldots ,0)^T$, $c^T\Sigma_\infty c = \sigma^2$, and we have
\begin{equation}\label{ar_n_upper}
\limsup_{\varepsilon\rightarrow 0}\varepsilon^2 \log E_{(x_0,\ldots , x_{n-1})}\tau_{(-1,1)} \le\frac{1}{2\sigma^2}.
\end{equation}
Thus, we have the desired upper bound for any starting points $x_0,\ldots ,x_{n-1}\in (-1,1)$. For the corresponding lower bound, we make the additional assumption that $x_0 =\ldots = x_{n-1}= 0$. Then $X_t^\varepsilon$ has a normal distribution with mean zero and variance
\[
\sigma_t^2 = \sum_{k=0}^{t-n}(B_{11}^k)^2 \le \sigma^2.
\]
Theorem \ref{thmlower} then immediately gives us
\begin{equation}\label{ar_n_lower}
\liminf_{\varepsilon\rightarrow 0}\varepsilon^2\log E_{(0,\ldots ,0)}\tau_{(-1,1)} \ge \frac{1}{2\sigma^2}.
\end{equation}
The upper and lower bounds in inequalities \ref{ar_n_upper} and \ref{ar_n_lower} together imply that
\begin{equation}
\lim_{\varepsilon\rightarrow 0}\varepsilon^2\log E_{(0,\ldots ,0)}\tau_{(-1,1)} = \frac{1}{2\sigma^2},
\end{equation}
and the proof is finished.

\medskip\noindent
{\bf Acknowledgements:} I would like to thank professor M. A. Lifshits for first showing me the method for lower bounds, and professor G\"{o}ran H\"{o}gn\"{a}s for many encouraging discussions. The financial support of the Academy of Finland (grant no. 127719) is gratefully acknowledged.


\end{document}